\documentclass[11pt]{amsart}
\usepackage[mathscr]{eucal}
\usepackage{amsfonts}
\usepackage{amsmath}
\usepackage{amsthm}
\usepackage{amssymb}
\usepackage{latexsym}

\theoremstyle{plain}
\newtheorem{definition}[equation]{Definition}

\newtheorem{lemma}[equation]{Lemma}
\newtheorem{proposition}[equation]{Proposition}
\newtheorem{theorem}[equation]{Theorem}
\newtheorem{example}[equation]{Example}

\theoremstyle{definition}
\newtheorem{remark}[equation]{Remark}

\numberwithin{equation}{subsection}

\begin{document}
\title{Koszulity and the Hilbert series of preprojective algebras}

\author{Pavel Etingof}
\address{Department of Mathematics, Massachusetts Institute of
Technology, Cambridge, MA 02139, U.S.A.}
\email{etingof@math.mit.edu}

\author{Ching-Hwa Eu}
\address{Department of Mathematics, Massachusetts Institute of
Technology, Cambridge, MA 02139, U.S.A.}
\email{ceu@math.mit.edu}

\maketitle{}

\section{Introduction}

The goal of this paper is to prove 
that if $Q$ is a connected non-Dynkin quiver then 
the preprojective algebra $\Pi_Q(k)$ of $Q$ over any field $k$ is
Koszul, and has Hilbert series $(1-Ct+t^2)^{-1}$, where 
$C$ is the adjacency matrix of the double $\bar Q$ of $Q$. 

We also prove a similar result for the partial preprojective algebra
$\Pi_{Q,J}(k)$ of any connected quiver $Q$, 
where $J\subset I$ is a subset of the set $I$ of
vertices of $Q$. By definition, $\Pi_{Q,J}(k)$ is the quotient of the path
algebra of $k\bar Q$ by the preprojective algebra relations
imposed only at vertices not contained in $J$. We show that 
if $J\ne \emptyset$ then $\Pi_{Q,J}(k)$ is Koszul, and its Hilbert series is 
$(1-Ct+D_Jt^2)^{-1}$, where $D_J$ is the diagonal matrix 
with $(D_J)_{ii}=0$ if $i\in J$ and $(D_J)_{ii}=1$, $i\notin J$. 

Moreover, we show that both results are valid in a slightly more
general framework of modified preprojective algebras, 
considered in \cite{K}. 

We note that the first result is known in most cases \cite{MV,MOV,O}.
In particular, it is known in general in characteristic zero
(\cite{MOV}), and in most positive characteristic cases \cite{MV,O}. 
Our argument, however, is elementary, and different from 
the arguments of \cite{MOV,O}, which are based on the theory of 
tensor categories. 

{\bf Acknowledgments.} P.E. is grateful to V. Ostrik for a useful discussion.
The work of P.E. was  partially supported by the NSF grant
 DMS-0504847 and by the CRDF grant RM1-2545-MO-03.

\section{Preliminaries} 

In this section we give some basic definitions and known results, 
which will be useful in the sequel. 

\subsection{Hilbert series of bimodules} Let 
$k$ be a field (of any 
characteristic). Let $I$ be a finite set, and 
let $R=\oplus_{i\in I}k$ be the algebra of $k$-valued
functions on $I$. 
An $R$-bimodule $W$ may be thought of as an $I\times I$-graded vector space 
$W=\oplus_{i,j\in I}W_{i,j}$. 
Then for two $R$-bimodules $U,\,W$ we have 
\begin{equation}\label{tenpro}
\left(W\underset{R}{\otimes}U\right)_{i,j}=\bigoplus
\limits_{s\in I}W_{i,s}\underset{k}{\otimes}U_{s,j}.
\end{equation}

The \emph{tensor algebra} $T_R(W)$ is defined as $\oplus_{m\ge
0}W^{\underset{R}{\otimes} m}$, with the tensor products over $R$.

Now let $W=\oplus_{d\geq0}W[d]$ be a $\mathbb Z_+$-graded
$R$-bimodule, with finite dimensional homogeneous subspaces. 
\begin{definition}
We define the \emph{Hilbert series} $h_W(t)$ to be a
matrix-valued series, with the entries
\begin{displaymath}
h_W(t)_{i,j}=\sum\limits_{d=0}^{\infty}dim\ W[d]_{i,j}t^d.
\end{displaymath}
\end{definition}

From (\ref{tenpro}), it follows that $h_{W\underset{R}{\otimes}U}(t)=h_W(t)h_U(t)$.

\subsection{Free products} 

Let $A,\,B$ be $k$-algebras containing $R$.

\begin{definition}
We define the \emph{free product} 
$A\underset{R}*B$ of $A$ and $B$ over $R$ to be the algebra
$T_R(A\oplus B)/J$, where $J$ is the ideal generated 
by elements 

1) $r_A-r_B$, where $r\in R$ and $r_A,r_B$ are its images in $A$
and $B$, respectively, and 

2) $a_1\otimes a_2-a_1a_2$, $a_1,a_2\in A$; $b_1\otimes
b_2-b_1b_2$, $b_1,b_2\in B$. 
\end{definition}

\begin{example}
If $A=T_R(V)/(r_i)$, $B=T_R(W)/(s_j)$ and $r_i, s_j$ are the
defining relations, then $A\underset{R}{*}B=T_R(V\oplus W)/(r_i,\,s_j)$.
\end{example}

\begin{example} Assume that we have decompositions of bimodules
$A=R\oplus A_+$, $B=R\oplus B_+$. Then 
\begin{equation}\label{astarb}
A\underset{R}{*}B=\bigoplus\limits_{m=0}^{\infty}
A\underset{R}\otimes\left(B_+\underset{R}\otimes
A_+\right)^m\underset{R}\otimes B
\end{equation}
\end{example}

\begin{lemma}\label{free product formula}
Suppose $A,B$ are $\mathbb Z_+$-graded algebras, with 
$A[0]=B[0]=R$. Let $h_A(t)=\frac{1}{1-\alpha}$,
$h_B(t)=\frac{1}{1-\beta}$. Then $h_{A\underset{R}{*}B}(t)=\frac{1}{1-\alpha-\beta}$.
\end{lemma}
\begin{proof} Let $A_+$, $B_+$ be the positive degree parts of
$A,B$. Then we have 
$$
h_{A_+}(t)=\frac{\alpha}{1-\alpha}, h_{B_+}(t)=\frac{\beta}{1-\beta},
$$
$$
h_{\oplus_{m\geq 0}\left(B_+\underset{R}\otimes
A_+\right)^m}(t)=\sum\limits_{m=0}^{\infty}\left(\frac{\beta}{1-\beta}\frac{\alpha}{1-\alpha}\right)^m=\left(1-\frac{\beta}{1-\beta}\frac{\alpha}{1-\alpha}\right)^{-1}.
$$
By (\ref{astarb}), we have
\begin{eqnarray*}
h_{A\underset{R}{*}B}(t)&=&(1-\alpha)^{-1}\left(1-\frac{\beta}{1-\beta}\frac{\alpha}{1-\alpha}\right)^{-1}(1-\beta)^{-1}\\
&=&\left((1-\beta)(1-\alpha)-\beta\alpha\right)^{-1}=\frac{1}{1-\alpha-\beta}.
\end{eqnarray*}
\end{proof}

\subsection{Quadratic and Koszul algebras} 

For the material in this subsection, we refer to \cite{PP}. 

Let $V$ be an $R$-bimodule, and $E\subset V\underset{R}\otimes V$
an $R$-subbimodule. To this data we may attach the 
\emph{quadratic algebra} $A:=T_R(V)/(E)$. 

\begin{definition} Assume that the bimodule $V$ is equipped with a
filtration by nonnegative integers. Then we define $A':=T_R({\rm
gr} V)/({\rm gr} E)$. 
\end{definition}

The following lemma is obvious. 

\begin{lemma}\label{aprimea}
We have the termwise inequality $h_{A'}(t)\geq h_A(t)$.
In other words, for every $i,j\in A$ and $d\ge 0$ we have 
$\dim A'[d]_{i,j}\ge \dim A[d]_{i,j}$.  
\end{lemma}

\begin{definition}
Let $A$ be a quadratic algebra with $A[0]=R$. 
$A$ is a \emph{Koszul algebra} if $Ext_A^i(R,R)$ 
(or dually $Tor_i^A(R,R)$) is concentrated in degree $i$.
\end{definition}
\begin{theorem}  \label{Hilbert inequality}
(The Golod-Shafarevich inequality)
Let $A$ be a $\mathbb{Z}_+$-graded algebra with $A[0]=R$, with
the bimodule of generators $V$ in degree $1$ and relations $E$ in
degree $2$. Let $C$ and $D$ be the matrices 
defined by $C_{i,j}=\dim V_{i,j}$, $D_{i,j}=\dim E_{i,j}$. Then:
\\
(i) If $\frac{1}{1-Ct+Dt^2}\geq 0$ termwise, then $h_A(t)\geq\frac{1}{1-Ct+Dt^2}$.\\
(ii) If $h_A(t)=\frac{1}{1-Ct+Dt^2}$, then $A$ is Koszul. 
\end{theorem}
\begin{proof}
From the exact Koszul complex, 
\begin{displaymath}
0\rightarrow K\rightarrow A\underset{R}{\otimes}E\rightarrow A\underset{R}{\otimes}V\rightarrow A\rightarrow R\rightarrow 0
\end{displaymath}
extended by the kernel $K$ on the left, using the
Euler-Poincar\'e principle, we obtain the equation
$$h_A(t)(1-Ct+Dt^2)-1=h_K(t)\geq 0,$$ so under the assumption of
(i) the inequality $h_A(t)-\frac{1}{1-Ct+Dt^2}\geq 0$ follows.
To prove (ii) note that in the above exact sequence,
$h_A(t)=\frac{1}{1-Ct+Dt^2}$ implies $h_K(t)=0$, therefore
$K=0$. Applying the functor 
 $R\underset{A}{\otimes}\_$ to this sequence (without the last
term $R$), we get the complex computing $Tor_i^A(R,R)$: 
\begin{displaymath}
0\overset{0}{\rightarrow}E\overset{0}{\rightarrow}V\overset{0}{\rightarrow}R\overset{0}{\rightarrow}0.
\end{displaymath}
So
\begin{displaymath}
Tor_i^A(R,R)=\left\{
\begin{array}{ll}
R&i=0\\
V&i=1\\
E&i=2\\
0&i>2,
\end{array}\right.
\end{displaymath}
and hence $A$ is Koszul. 
\end{proof}

\section{Hilbert series of preprojective algebras} 

\subsection{Partial preprojective algebras}

Let $Q$ be a finite quiver with vertex set $I$
(we allow multiple edges and self-loops).
For each $a\in Q$ let $h(a),t(a)\in I$ denote the head and tail 
of the edge $a$. Let $V$ be the $k$-vector space spanned by the edges of $Q$; 
it is naturally an $R$-bimodule, where $R=\oplus_{i\in I}k$. 
Recall that the path algebra $kQ$ is the algebra $T_RV$. 
 
Let $\bar Q$ the \emph{double} of $Q$, obtained from 
$Q$ by keeping the same vertex set and adding
a new edge $a^*$ for $j$ to $i$ for each edge $a$ from $i$ to
$j$. 

Let $J\subset I$ be a subset of ``white'' vertices (the other
vertices are colored black).

\begin{definition}
We define the \emph{partial preprojective algebra}
\begin{displaymath}
\Pi_{Q,J}(k)=k{\bar Q}/\left(\sum\limits_{a\in Q, h(a)\in I\setminus J}aa^*-\sum\limits_{a\in Q, t(a)\in I\setminus J}a^*a\right)
\end{displaymath}
In particular, if $J=\emptyset$ then we write 
$\Pi_{Q}(k)$ for $\Pi_{Q,\emptyset}(k)$ and call it the
\emph{preprojective algebra} of $Q$. 
\end{definition}

\begin{remark} In a similar way, for any commutative ring 
$K$ one defines the algebra $\Pi_{Q,J}(K)$. When no confusion is
possible, we will denote $\Pi_{Q,J}(K)$ by $\Pi_{Q,J}$.  
\end{remark} 

\subsection{The Hilbert series for an extended Dynkin quiver} 

Let $Q$ be an extended Dynkin quiver (in particular, we allow 
the quiver $\widetilde{A}_0$, which has one vertex and one edge). 

\begin{proposition}\label{prop1}
The Hilbert series $h_{\Pi_Q(k)}(t)$ equals $\frac{1}{1-Ct+t^2}$, where $C$ is the adjacency matrix of $\bar Q$.
\end{proposition}
\begin{proof}
First we show this for char$k=0$. (In this case the result is well
known, and the argument actually works if $p$ is a good prime for
$Q$, in particular for any $p>5$). 

Let $\Gamma\subset SL_2(k)$ be the finite group attached to $Q$
via the McKay's correspondence, $f_i$ the primitive 
idempotents of the irreducible representations of $\Gamma$, 
$f=\sum\limits_if_i$, and $A=\Gamma\ltimes k[x,y]$ the skew group algebra.

Then we have the following result, due to G. Lusztig \cite{L}
(see also \cite{CBH}). 

\begin{proposition}\label{cbh} 
\begin{displaymath}
\Pi_Q(k)=\bigoplus_{i,j}f_i A f_j=fAf
\end{displaymath}
\end{proposition}

By applying the functor $\_\otimes k\Gamma$ to the exact Koszul complex of $k[x,y]$ (over $k$) we obtain the exact complex $\tilde{K}^{\bullet}$:
\begin{displaymath}
0\longrightarrow A(2)\longrightarrow A\underset{k\Gamma}
{\otimes} A[1]\longrightarrow A\longrightarrow k\Gamma\longrightarrow 0
\end{displaymath}
(here $(2)$ denotes the shift in grading by $2$). 
By Proposition \ref{cbh}, 
the complex $K^{\bullet}=f\tilde{K}^{\bullet}f$ has the form
\begin{displaymath}
0\longrightarrow \Pi_Q(k)(2)\longrightarrow \Pi_Q(k)\underset{R}{\otimes} {\Pi_Q(k)[1]}\longrightarrow\Pi_Q(k)\longrightarrow R\longrightarrow 0.
\end{displaymath}
We conclude that this complex is exact, and from the
Euler-Poincar\'e principle obtain the equation 
$$1=h_R(t)=h_{\Pi_Q(k)}(t)-h_{\Pi_Q(k)}(t)Ct+h_{\Pi_Q(k)}(t)t^2.$$

Now consider the case char$k\neq 0$.
Let $T$ be the torsion part of $\Pi_Q(\mathbb{Z})$. 
Then $\Pi_Q(\mathbb{Z})/T$ is a free $\mathbb{Z}$-module, 
and by the characteristic zero result has the Hilbert series
$\frac{1}{1-Ct+t^2}$ (indeed, it suffices to take tensor product
with $\mathbb C$). 
Further, $T\underset{\mathbb{Z}}{\otimes} k\subset\Pi_Q(k)$ is an ideal, and
the quotient algebra \linebreak
$\Pi_Q(k)/T\underset{\mathbb{Z}}{\otimes} k=
\left(\Pi_Q(\mathbb{Z})/T\right)\underset{\mathbb{Z}}{\otimes} k$
has the same Hilbert series. Since $C$ is an extended Cartan
matrix, its largest eigenvalue is $2$ and hence we get that
$\Pi_Q(k)/T\underset{\mathbb{Z}}\otimes k$ has Gelfand-Kirillov dimension
$2$. By~\cite{BGL}, $\Pi_Q(k)$ is a prime Noetherian
algebra of Gelfand-Kirillov dimension $2$.
But for any prime Noetherian algebra $A$ of Gelfand-Kirillov
dimension $d$ and any nonzero two-sided ideal $I\subset A$, one
has GKdim$(A/I)\leq d-1$,~\cite[Corollary 8.3.6.]{MR}. Therefore,
we see that $T\underset{\mathbb{Z}}{\otimes} k=0$, and hence $T=0$, as
desired. 
\end{proof}

\subsection{Star-shaped quivers}

In this subsection we prove the Hilbert series formula 
for star-shaped quivers with node being a white vertex.  

\begin{lemma}\label{stars}
Let $Q$ be a quiver 
with vertex set $I=\lbrace{1,...,n+1\rbrace}$ and arrows
$n+1\xrightarrow{a_{i,1},\ldots a_{i,r_i}}i\in I$, 
and with a set of white vertices $J=\{i_1,\ldots i_m,n+1\}$. 
Then $h_{\Pi_{Q,J}}(t)=\frac{1}{1-Ct+D_Jt^2}$ where $D_J$ 
is a diagonal matrix with $(D_J)_{ii}=\left\{
\begin{array}{cc}1&i\notin J\\
0&i\in J.\end{array}\right.$
\end{lemma}
\begin{proof}
First, we observe that $\Pi_{Q,J}=A_1*\ldots*A_{n+1}$ where $A_i$
is the partial preprojective algebra of the quiver $Q_i$ with
vertex set $I$ and only the arrows
$n+1\xrightarrow{a_{i,1},\ldots a_{i,r_i}} i$ (if $i=n+1$, these
arrows are self-loops). By Lemma
\ref{free product formula}, it is enough to prove the result 
in the case when this free product has only one factor, i.e. 
for $n=1$, and for $n=2$, $r_2=0$. 

If all vertices are white, then $\Pi_{Q,J}$ is a path
algebra and the result is clear from counting paths. 
So it remains to consider the
case $n=2$, $r_2=0$, where $1$ is a black vertex. 
In this case, the quiver $Q$ has edges $a_1,...,a_r$ going from
$2$ to $1$, and we have one defining relation
$$
a_1a_1^*+...+a_ra_r^*=0. 
$$
Denote the algebra $\Pi_{Q,J}$ in this case by $A(r)$. 

Consider first the case $r=1$. In this case 
the algebra $A(r)$ has only one quadratic element 
$a_1^*a_1$ (up to scaling) and no cubic elements, so the formula
easily follows. 

Now consider the case $r>1$. Let us introduce the filtration on $A:=A(r)$ by setting 
$\deg(a_1)=\deg(a_1^*)=1$, and $\deg(a_i)=\deg(a_i^*)=0$ for
$i>1$. In this case $A'=A(1)\underset{R}{*}B$, where $B$ is the
path algebra of the quiver with edges $a_i,a_i^*$, $i=2,...,r$.
It then follows from Lemma \ref{free product formula} and the
$r=1$ case that the desired Hilbert series formula holds for $A'$. 
By Lemma \ref{aprimea}, this implies that 
$$
h_A(t)\le \frac{1}{1-Ct+D_Jt^2}.
$$
Then by Theorem \ref{Hilbert inequality}, we see that 
$$
h_A(t)=\frac{1}{1-Ct+D_Jt^2}, 
$$
as desired. 
\end{proof}

\subsection{Main results}

\begin{theorem}\label{parprep}
Let $Q$ be a connected quiver and a nonempty set of 
white vertices. Then $h_{\Pi_{Q,J}}(t)=\frac{1}{1-Ct+D_Jt^2}$.
In particular, $\Pi_{Q,J}$ is a Koszul algebra. 
\end{theorem}

\begin{proof}
The second statement follows from the first statement and Theorem \ref{Hilbert
inequality}. So we only need to prove the first statement. 

We will prove by induction in the number of vertices  
that the statement is true for any quiver $Q$ whose 
every connected component contains a white vertex.

If $Q$ has just one vertex, the statement is clear.

To make the induction 
step, assume the formula is right for $\leq n$ vertices, and that
$Q$ has $n+1$ vertices. 

Let $I$ be the vertex set of $Q$. 
Select any white vertex $w$ in $I$ and consider the
subquiver $Q_0$ with vertex set $I_0\subset I$, consisting of $w$
and the vertices adjacent to it, whose arrows are all the arrows
of $Q$ which touch $w$. Also, let
$Q'$ be the quiver with the vertex set $I\setminus\{w\}$, where
the vertices in $I_0$ are colored white and the other ones are
colored in the same way as in $Q$, and with the set of arrows $Q\setminus Q_0$. 
Finally, let $\widehat{Q}_0$, $\widehat{Q}'$ be the quivers with
vertex set $I$ and with the same arrows as in $Q_0,Q'$
respectively, and let $J_0,J'$ be the sets of white vertices of
$Q_0$ and $Q'$. 

Introduce a filtration on $\Pi_{Q,J}$, by setting a grading, 
such that the arrows inside $\bar Q_0$ have degree $1$ 
and the other ones have degree $0$. Then 
$\Pi_{Q,J}'=\Pi_{\widehat{Q}_0,J_0}*\Pi_{\widehat{Q}',J'}$.
Therefore, by Lemma \ref{aprimea}, 
we get 
$$
h_{\Pi_{Q,J}}(t)\leq
h_{\Pi_{Q,J}'}(t)=h_{\Pi_{\widehat{Q}_0,J_0}*\Pi_{\widehat{Q}',J'}}(t).
$$
Now, by Lemma \ref{stars} for $Q_0$, 
$h_{\Pi_{\widehat{Q_0},J_0}}(t)=\frac{1}{1-C_0t+D_0t^2}$ where $C_0$
is the adjacency 
matrix of $\bar Q_0$, and $D_0$ is the diagonal matrix, such that
$(D_0)_{ii}=1$ 
if $i$ is a black vertex in $I_0$, and $0$ otherwise.
Also, applying the induction assumption to $Q'$, we find 
$h_{\Pi_{\widehat{Q'},J'}}(t)=\frac{1}{1-(C-C_0)t+D't^2}$, where
$D'$ is the diagonal matrix with $(D')_{ii}=1$ if $i$ is a black
vertex $\notin I_0$, and $0$ otherwise. 
Therefore, by Lemma \ref{free product formula} 
we have $h_{\Pi_{Q,J}}(t)\leq\frac{1}{1-Ct+D_Jt^2}$. 
From this, it follows $\frac{1}{1-Ct+D_Jt^2}\geq 0$. Now, by
Theorem \ref{Hilbert inequality}, the result follows. 
\end{proof}

\begin{theorem}\label{nondyn}
The Hilbert series 
$h_{\Pi_Q}(t)$ 
for a connected non-Dynkin quiver $Q$ is $\frac{1}{1-Ct+t^2}$. In
particular, 
$\Pi_Q$ is Koszul.
\end{theorem}

\begin{proof}
Again, the second statement follows from the first one and Theorem \ref{Hilbert
inequality}, and we only prove the first statement. 

For this, we use the following easy (and well known) lemma. 

\begin{lemma}
Any connected non-Dynkin quiver $Q$ contains an extended Dynkin subquiver $Q_{E}$.
\end{lemma}

Let $I_{E}$ be the vertex set of $Q_{E}$ (it is possible that 
an arrow $a$ between $i,j\in I_{E}$ belongs to $Q$ but not to
$Q_{E}$). Let $Q'$ be the quiver with
vertex set $I$, such that the vertices in $I_{E}$ are white and
the other ones are black, and the set of arrows $Q\setminus
Q_{E}$. Then every connected component of $Q'$ contains at least one
white vertex. Introduce a filtration on $\Pi_Q$ by setting a grading, such that the arrows in
$\bar Q_{E}$ have degree $1$ and the other ones have degree $0$. Then
$\Pi_Q'=\Pi_{\widehat{Q}_{E}}*\Pi_{Q',I_E}$, where 
$\widehat{Q}_{E}$ is the quiver $Q_{E}$ with adjoined vertices 
of $I\setminus I_{E}$. 

Now, by Proposition \ref{prop1}, 
$h_{\Pi_{\widehat{Q}_{E}}}=\frac{1}{1-{C_{E}}t+D_{E}t^2}$
where $D_{E}$ is the diagonal matrix, such that
$\left(D_{E}\right)_{ii}$ is $1$ if $i\in I_{E}$ and $0$
otherwise, and $C_E$ is the adjacency matrix of 
of the double of $\widehat{Q}_E$. 
Also, by Theorem \ref{parprep},
$h_{\Pi_{Q',I_E}}=\frac{1}{1-(C-C_{E})t+(1-D_{E})t^2}$. Hence, 
by Lemma \ref{free product formula}, we obtain 
$$
h_{\Pi_Q}\leq
h_{\Pi_Q'}=h_{\Pi_{\widehat
{Q}_{E}}*\Pi_{Q',I_E}}=
\frac{1}{1-Ct+t^2}.
$$ 
From this, it follows that $\frac{1}{1-Ct+t^2}\geq 0$. 
Hence, Theorem \ref{Hilbert inequality} implies the result. 
\end{proof}

\subsection{Modified preprojective algebras} 

It turns out that our results hold for a slightly 
more general class of preprojective algebras. 
Namely, let $Q$ be a quiver with vertex set $I$, and $\bar Q$ its
double. Let $J\subset I$. Let $\gamma$ be a $k^\times$-valued function on the set of edges 
of $\bar Q$ which begin or end at $I\setminus J$. Define the 
modified partial preprojective algebra
$\Pi_{Q,J}^\gamma(k)$ to be the quotient 
of $k\bar Q$ by the relation 
$$
\sum\limits_{a\in Q, h(a)\in I\setminus J}\gamma_a aa^*-\sum\limits_{a\in
Q, t(a)\in I\setminus J}\gamma_{a^*} a^*a=0.
$$
Obviously, this is a generalization of the usual partial preprojective algebras
$\Pi_{Q,J}(k)$, which are obtained if $\gamma=1$. 

If $J=\emptyset$, $\Pi_{Q,J}^\gamma(k)$ 
is called a modified preprojective algebra 
and denoted by $\Pi_Q^\gamma(k)$. Such algebras are considered in
\cite{K}. 

\begin{theorem}\label{gamma} Theorem \ref{parprep}
and Theorem \ref{nondyn} hold for the algebras
$\Pi_{Q,J}^\gamma(k)$
and $\Pi_Q^\gamma(k)$. 
\end{theorem}

\begin{proof} The proof of Theorem \ref{parprep} carries out
verbatim to the case of general $\gamma$. 

To prove Theorem \ref{nondyn} for $\Pi_Q^\gamma$, 
we only need to consider the extended Dynkin case, since it is
the only case when the proof needs to be changed. 
Also, note that if $Q$ is a tree, then 
the algebras $\Pi_Q^\gamma$ are pairwise isomorphic for all
choices of $\gamma$ (by rescaling the edges). So it is necessary
to consider only the case of type $\widetilde {A}_{n-1}$.

In this case the edges are $a_1,...,a_n$ and $a_1^*,...,a_n^*$ 
(where $a_i$ goes from $i$ to $i+1$), and the defining relations are
$$
\gamma_{a_i^*}a_i^*a_i-\gamma_{a_{i-1}}a_{i-1}a_{i-1}^*=0,
$$
where $i-1$ is computed modulo $n$. 

These relations show that $\Pi_Q^\gamma$ is spanned by paths 
in which there is no expressions $a_i^*a_i$. Counting such paths,
we easily get that 
$$
h_{\Pi_Q^\gamma}\le (1-Ct+t^2)^{-1},
$$ 
which by Theorem \ref{Hilbert inequality} implies the result. 
\end{proof}


\begin{thebibliography}{999}

\bibitem[BGL]{BGL} D. Baer, W. Geigle, and H. Lenzing, The preprojective algebra of
a tame hereditary algebra, Comm. Algebra 15 (1987), 425-457.

\bibitem[CBH]{CBH} 
W. Crawley-Boevey and M. P. Holland: \emph{Noncommutative deformations of Kleinian singularities}, Duke Mathematical Journal, Vol. 92, No. 3 (1998)

\bibitem[K]{K} M. Kleiner, 
The graded preprojective algebra of a quiver, 
Bulletin of London Math. Soc., v.36 (2004), p.13-22. 

\bibitem[L]{L}
G. Lusztig,  Quivers, perverse sheaves, 
and quantized enveloping algebras, J. Amer. Math. Soc. 4 (1991),
365--421. 

\bibitem[MOV]{MOV} A. Malkin, V. Ostrik, M. Vybornov,
Quiver varieties and Lusztig's algebra,
math.RT/0403222.

\bibitem[MR]{MR}
J.C. McConnell and J.C. Robson: \emph{Noncommutative Noetherian Rings}, Graduate Studies in Mathematics, Vol. 30 (2001)

\bibitem[MV]{MV}
R. Martinez-Villa, Applications of Koszul algebras: the preprojective algebra,
in Representation theory of algebras, Cocoyoc, 1994, CMS
Conf. Proceedings, v. 18, AMS, Providence, RI pp. 487-504.

\bibitem[O]{O} V. Ostrik, 
Module categories over representations of $SL_q(2)$
and graphs II, math.QA/0509530. 

\bibitem[PP]{PP} A. Polishchuk, L. Positselski, Quadratic
algebras, AMS, 2005. 
\end{thebibliography}
\end{document}